\documentclass[11pt]{article}
\usepackage{enumerate}
\usepackage{amsmath}
\usepackage{amssymb}
\usepackage{amsthm}
\usepackage{amscd}
\DeclareSymbolFontAlphabet{\Bbb}{AMSb}
\newcommand{\id}{\text{{\rm id}}}
\newcommand{\Span}{\text{{\rm span}}}

\newcommand{\norm}[1]{ \|#1 \| }
\newcommand{\K}{\Bbb{K}}
\newcommand{\C}{\Bbb{C}}

\newcommand{\N}{\Bbb{N}}

\newcommand{\LL}{{\cal L}}
\newcommand{\MM}{{\cal M}}
\newcommand{\PP}{\Pi}

\newcommand{\On}{\mathcal{O}}
\newcommand{\type}{\mathbf{T_2}}
\newcommand{\cotype}{\mathbf{C_2}}
\newcommand{\con}{\mathbf{M^{(2)}}}
\newcommand{\conu}{\mathbf{M^{(u)}}}

\newcommand{\concave}{\mathbf{M_{(2)}}}

\newcommand{\Id}{\hookrightarrow}

\newcommand{\ui}{\mathcal{S}}
\theoremstyle{definition}

\newtheorem{rem}{Remark}

\theoremstyle{plain}
\newtheorem{lemma}{Lemma}
\newtheorem{theo}{Theorem}
\newtheorem{cor}{Corollary}

\theoremstyle{remark}
\numberwithin{equation}{section}
\begin{document}
\bibliographystyle{amsalpha}
\title{\bf Bennett--Carl Inequalities for Symmetric Banach Sequence Spaces and
 Unitary Ideals}
\author{Andreas Defant and Carsten Michels 
\thanks{The second named author was supported by a GradF\"oG Stipendium
 of the Land Niedersachsen (State of Lower Saxony)} \\
 \small \sl Carl von Ossietzky Universit\"at Oldenburg,
 Fachbereich Mathematik \\
 \small \sl Postfach 2503, D-26111 Oldenburg, Germany \\
  \small \tt defant\rm @\tt mathematik.uni-oldenburg.de \\
  \small \tt    michels\rm @\tt mathematik.uni-oldenburg.de }
\date{}
\maketitle
\begin{abstract}
We prove an abstract interpolation theorem which interpolates
 the $(r,2)$-summing and $(s,2)$-mixing norm of a fixed operator in the
 image and range space. Combined with interpolation formulas for spaces
 of operators 
 we  obtain as an application the original Bennett--Carl inequalities for 
identities acting between Minkowski spaces $\ell_u$
 as well as their analogues for 
 Schatten classes $\ui_u$. Furthermore, our techniques motivate 
a  study of Bennett--Carl inequalities
within a more general setting of symmetric Banach sequence spaces and 
unitary ideals.
\end{abstract}
In \cite{bennett} and \cite{carl} Bennett and Carl independently
proved
 the following inequalities: For $1\le u \le 2$ and $1 \le u \le v \le 
\infty$ the identity operator $\id : \ell_u \Id \ell_v$ is absolutely
$(r,2)$-summing, i.\,e. there is a constant $c>0$ such that for each set
 of finitely many $x_1, \dots ,x_n \in \ell_u$
$$ \bigl ( \sum_{k=1}^n \norm{x_k}_{\ell_v}^r \bigr )^{1/r} \le
 c \cdot \sup_{\norm{x'}_{\ell_{u'}}\le 1} \bigl ( \sum_{k=1}^n | \langle
  x',x_k \rangle |^2 \bigr )^{1/2},$$
if (and only if) $1/r \le 1/u - \max(1/v, 1/2)$. 
\\[10pt]
This result improved upon older ones of Littlewood and Orlicz, and is nowadays 
of extraordinary importance in the 
theory of eigenvalue distribution of power compact operators (see e.\,g.
 \cite{koenig} and \cite{pietsch2}). Later in \cite{CD92} the 
``Bennett--Carl inequalities'' were extended within the setting of so-called
 mixing operators (originally invented by Maurey \cite{maurey}): For
 $1\le u \le 2$ and $1 \le u \le v \le \infty$ every $s$-summing
operator $T$ defined on $\ell_v$ has a $2$-summing restriction
 to $\ell_u$ if (and only if) $1/s \ge 1/2 - 1/u + \max(1/v,1/2)$.
\\[5pt]
The crucial step in the proofs of Bennett and Carl is to establish the case 
$1 \le u \le v =2$ which in terms of finite-dimensional spaces reads 
as follows: 
For $1 \le u \le v= 2$
 and $2 \le r \le \infty$
 such that  $1/r = 1/u - 1/2$
\begin{equation}
\label{rlulv}
\sup_{n} \pi_{r,2}(\ell_u^n \Id \ell_2^n) < \infty,
\end{equation}
where $\pi_{r,2}(\ell_u^n \Id \ell_2^n)$ denotes the $(r,2)$-summing norm of the
 embedding $\ell_u^n \Id \ell_2^n$. 
\par
Note that the formula $1/r=1/u-1/2$
 occurs by a ``naive interpolation'' of the parameter $r$ 
between the two well-known border cases
\begin{align*}
\sup_n \pi_{\infty,2}(\ell_2^n \Id \ell_2^n)  = \sup_n \norm{\ell_2^n 
\Id \ell_2^n} < \infty \\
\sup_n \pi_{2,2}(\ell_1^n \Id \ell_2^n)  = \sup_n \pi_2(\ell_1^n \Id 
\ell_2^n) < \infty
\end{align*}
($\pi_2$ the $2$-summing norm): For $1 \le u \le 2$ choose $0 \le \theta 
\le 1$ with 
$1/u =(1-\theta)/2 + \theta/1$, then $1/r = (1-\theta)/\infty + \theta/2 =
 1/u-1/2$.
\\[10pt]
Nevertheless the literature so far does not offer an  approach
 to the Bennett--Carl inequalities
 within the framework of interpolation theory. We prove an abstract 
interpolation formula for the mixing norm of a fixed operator, and
 obtain as an application not only the original 
Bennett--Carl inequalities but also their ``non-commutative'' analogues for 
finite-dimensional Schatten classes. Moreover, we consider 
Bennett--Carl inequalities in a more general setting of symmetric Banach
 sequence spaces and unitary ideals, and apply these results to Orlicz and
 Lorentz sequence spaces. 
\par
 For further extensions of the Bennett--Carl inequalities within the setting
 of  Orlicz sequence spaces see a recent paper of Maligranda 
and Masty\l o \cite{mm99}.   
\section{\sc Preliminaries}
If $(a_n)$ and $(b_n)$ are scalar sequences we write \mbox{$a_n \prec b_n$}
 whenever there is some $c \ge 0$ such that $a_n \le c \cdot b_n$ for all $n$, 
  and 
\mbox{$a_n \asymp b_n$} whenever \mbox{$a_n \prec b_n$} and 
\mbox{$b_n \prec a_n$}.
 For $1 \le p \le \infty$ the number $p'$ is defined by $1/p + 1/{p'} =1$.
\par
 We shall use standard notation and
 notions from Banach space theory, as presented e.\,g. in 
\cite{djt}, \cite{lt} and \cite{tj}.
If $E$ is a Banach space, then
 $B_E$ is its (closed) unit ball and $E'$ its dual. We consider complex
 Banach spaces only (but note that most of the consequences of our abstract
interpolation results can be formulated also for real spaces).
As usual $\LL(E,F)$
 denotes the Banach space of all (bounded and linear) operators from
 $E$ into $F$ endowed with the operator norm $\norm{\cdot}$.
With $\type(E)$ and $\cotype(E)$  we denote the (Gaussian) 
type~2 
and cotype~2  constant of a Banach space $E$ with this
property, respectively, and with
 $\mathbf{M_{(u)}}(E)$ and $\conu(E)$ the $u$-concavity and $u$-convexity 
constant of a Banach lattice $E$.
\par 
We call a Banach space $E \subset c_0$ (the space of all zero sequences) 
a symmetric Banach sequence space if the $i$-th standard unit vectors $e_i$ 
form a symmetric basis, i.\,e. the $e_i$'s form a Schauder basis such that 
$\norm{x}_E = \norm{\sum_{i=1}^\infty \varepsilon_i x_{\pi(i)} e_i}_E$
for each $x \in E$, each permutation $\pi$ of $\N$ and each choice of 
scalars $\varepsilon_i$ with $|\varepsilon_i|=1$. Moreover, denote 
for each $n$ the subspace 
$\Span \{e_i \, | \, 1 \le i \le n\}$ of $E$ by $E_n$.  
 Together 
with its natural order a symmetric Banach sequence space $E$ forms a  
Banach lattice, and clearly its basis is $1$-unconditional. 
The associated unitary
 ideal $\ui_E$ is the Banach space of all compact operators $T\in \LL(\ell_2
,\ell_2)$ with singular numbers $(s_i(T))_i$ 
 in $E$ endowed with the norm $\norm{T}_{\ui_E}:= 
\norm{(s_i(T))_i}_E$; with $\ui_E^n$ we denote $\LL(\ell_2^n,\ell_2^n)$
 together with the norm $\norm{T}_{\ui_E^n}:= 
\norm{(s_i(T))_{i=1}^n}_{E_n}$. For $E=\ell_u$ ($1 \le u < \infty$) one gets
 the well-known Schatten-$u$-class $\ui_u$; for simplicity put 
$\ui_\infty:= \LL(\ell_2,\ell_2)$.
\par
For all information on Banach operator ideals, in particular on summing
 and mixing operators, see e.\,g. \cite{df}, \cite{djt} and \cite{pietsch}.
 An operator $T \in \LL(E,F)$ is called absolutely $(r,p)$-summing 
$(1 \le p \le r \le \infty)$ if there is a constant $\rho \ge 0$ such that
$$
\bigl ( \sum_{i=1}^n \norm{Tx_i}^r \bigr )^{1/r} \le \rho
\cdot \sup \biggl \{ \bigl ( \sum_{i=1}^n |\langle x',x_i \rangle |^p 
\bigr )^{1/p} \, | \, x' \in B_{E'} \biggr \}
$$
 for all finite sets of elements $x_1, \dots, x_n \in E$
 (with the obvious modifications for $p$ or $r$ $=\infty$). In this case,
 the infimum over all possible $\rho \ge 0$ is denoted by $\pi_{r,p}(T)$, 
and the Banach operator ideal of all absolutely $(r,p)$-summing 
operators by $(\PP_{r,p},\pi_{r,p})$; the special case $r=p$ gives
 the ideal $(\PP_p, \pi_p)$ of all absolutely $p$-summing operators.
\par
An operator $T \in \LL(E,F)$ is called $(s,p)$-mixing 
$(1 \le p \le s \le \infty)$ whenever its composition with an 
arbitrary operator $S \in \PP_s(F,Y)$ is absolutely $p$-summing; with
 the norm
$$ \mu_{s,p}(T) := \sup\{\pi_p(ST) \, | \, \pi_s(S) \le 1 \} $$
the class $\MM_{s,p}$ of all $(s,p)$-mixing operators forms again a 
 Banach operator ideal. Obviously, $(\MM_{p,p}, \mu_{p,p}) =
 (\LL, \norm{\cdot})$ and $(\MM_{\infty,p}, \mu_{\infty,p}) = (\PP_p, \pi_p)$.  
 Recall that due to \cite{maurey} (see also \cite[32.10--11]{df}) summing 
and mixing operators are closely related: 
$$ (\MM_{s,p}, \mu_{s,p}) \subset (\PP_{r,p}, \pi_{r,p}) \qquad \text{for } 
   1/s + 1/r = 1/p, $$
and ``conversely'' 
$$  (\PP_{r,p}, \pi_{r,p}) \subset (\MM_{s_0,p}, \mu_{s_0,p}) \qquad
 \text{for } 1 \le p \le s_0 < s \le \infty \text{ and } 1/s + 1/r =
1/p.  $$
Moreover, it is known that each 
$(s,2)$-mixing operator on a cotype~2 space is even $(s,1)$-mixing
 (see again \cite{maurey} and \cite[32.2]{df}).
\par
For an operator $T \in \LL(E,F)$ 
 the $n$-th Weyl number $x_n(T)$ of $T$ is defined by 
$$ x_n(T) := \sup \{a_n(TS) \,| \, S \in \LL(\ell_2, E) \text{ with }
  \norm{S}=1 \},$$
where $a_n(TS)$ denotes the $n$-th approximation number of $TS$.
We will use the following important inequality of K\"onig in order to 
obtain lower estimates:
\begin{equation}
\label{weyl}
n^{1/r} \cdot x_n(T) \le \pi_{r,2}(T), \qquad T \in \PP_{r,2}
\end{equation}
(for all details on $s$-numbers and this inequality see \cite[2.a.3]{koenig} 
or \cite{pietsch2}).
\section{\sc Complex Interpolation of Mixing Operators}
The aim of this section is to prove a complex interpolation formula for the 
 mixing norm of a fixed operator acting between two complex interpolation
 spaces.
\par
For all information on complex interpolation  we refer to \cite{BL}.
A couple $[E_0,E_1]$ of Banach spaces is called an interpolation couple
 if there exists a topological Hausdorff vector space $E$ in which 
$E_0$ and $E_1$ can be continuously embedded. 
 We speak of a finite-dimensional interpolation couple $[E_0, E_1]$, if
 $E_0$ and $E_1$ are finite-dimensional Banach
 spaces with the same dimensions. 
 For an interpolation couple $[E_0,E_1]$ of Banach
 spaces and $0 \le \theta\le 1$ we denote by 
$[E_0,E_1]_\theta = E_\theta$ the interpolation space 
obtained by the complex interpolation method of Calder\'on. 
Well-known examples 
of complex interpolation spaces are the Minkowski spaces $\ell_p^n(E)$ and 
Schatten classes $\ui_p^n$: For $1\le p_0,p_1 \le \infty $, 
$0 < \theta < 1$ and an
 interpolation couple $[E_0,E_1]$
$$ [\ell_{p_0}^n(E_0), \ell_{p_1}^n(E_1)]_\theta = \ell_p^n(E_\theta)
\qquad \text{and} \qquad [\ui_{p_0}^n,\ui_{p_1}^n]_\theta = \ui_p^n$$
(isometrically), where $1/p = (1-\theta)/{p_0} + \theta/{p_1}$ (for the complex 
interpolation formula for $\ell_p^n(E)$'s see
 \cite[5.1.2]{BL}, whereas the formula for Schatten 
classes can be deduced from e.\,g. \cite[Satz~8]{pt} and 
the complex reiteration theorem \cite[4.6.1]{BL}). 
For $0 \le \theta <1$  
a $\theta$-Hilbert space is an interpolation space $[E_0,E_1]_\theta$ 
where $E_1$ is a Hilbert space (this notion goes back to Pisier); 
in particular,  $\ell_p^n$ and $\ui_p^n$ for $1 < p < \infty$
are $\theta$-Hilbert spaces for $\theta=1-|1-{2/p}|$.
\par
The following complex interpolation theorem for the mixing norm is our main
 abstract tool.
\begin{theo}
\label{mixing}
Let $2 \le s_0,s_1 \le \infty$, $0 \le \theta \le 1$ and $s_\theta$ 
given by $1/{s_\theta}=(1-\theta)/{s_0} + \theta / {s_1}$. Then for two
 finite-dimensional interpolation couples $[E_0,E_1]$, $[F_0,F_1]$
 and each $T \in \LL(E_\theta,F_\theta)$
$$ \mu_{s_\theta,2}(T : E_\theta \rightarrow F_\theta) \le d_\theta[E_0,E_1]
 \cdot \mu_{s_0,2}(T : E_0 \rightarrow F_0)^{1-\theta} \cdot
\mu_{s_1,2}(T : E_1 \rightarrow F_1)^{\theta},$$
where
\begin{equation}
 d_\theta[E_0,E_1]:= \sup_m \norm{\LL(\ell_2^m, E_\theta) 
\hookrightarrow [\LL(\ell_2^m, E_0), \LL(\ell_2^m, E_1)]_\theta}. 
\end{equation}
\end{theo}
\proof Consider for $\eta=0,\theta,1$ the bilinear mapping
$$
\begin{array}{lccccc}
\Phi_\eta^{n,m}: &\ell_{s_\eta}^n(F_\eta') & \times 
&\LL(\ell_2^m,
E_\eta) &\longrightarrow & \ell_2^m(\ell_{s_\eta}^n) \\
&(y_1', \dots , y_n') &\times & S &\longmapsto
& ((\langle y_k',TSe_j \rangle )_k)_j
\end{array},
$$
where $(e_j)$ denotes the canonical basis in $\C^m$. By the discrete 
characterization of the mixing norm (see \cite{maurey} or 
\cite[32.4]{df}) $\mu_{s_\eta}(T: E_\eta \rightarrow F_\eta)$ is the 
infimum over all $c \ge 0$ such that for all $n,m$, all $y_1', \ldots,
 y_n' \in F_\eta'$ and all $x_1, \ldots, x_m \in E_\eta$
$$
\left ( \sum_{j=1}^m \left ( \sum_{k=1}^n |\langle y_k',Tx_j \rangle |
^{s_\eta} \right )^{2/{s_\eta}} \right )^{1/2} \le 
c \cdot  \left ( \sum_{k=1}^n \norm{y_k'}_{F_\eta'}
^{s_\eta} \right )^{1/{s_\eta}} \cdot
\sup_{x' \in B_{E_\eta'}} \left ( \sum_{j=1}^m |\langle x',x_j 
\rangle |^2 \right )^{1/2}.
$$
Since for each $S = \sum_{j=1}^m e_j \otimes x_j \in \LL(\ell_2^m,E_\eta)$
$$
\norm{S} = \sup_{x' \in B_{E_\eta'}} \left ( \sum_{j=1}^m | \langle x',x_j
 \rangle |^2 \right )^{1/2},
$$
we obviously get that
$$
\mu_{s_\eta,2}(T:E_\eta \rightarrow F_\eta) = \sup_{n,m} 
\norm{\Phi_\eta^{n,m}}.
$$
Now the proof follows by bilinear complex interpolation: For the interpolated
 bilinear mapping
$$ [\Phi_0^{n,m},\Phi_1^{n,m}]_\theta : [\ell_{s_0}^n(F_0'), \ell_{s_1}^n(F_1')]_\theta 
   \times [\LL(\ell_2^m,E_0), \LL(\ell_2^m,E_1)]_\theta \longrightarrow
   [\ell_2^m(\ell_{s_0}^n),\ell_2^m(\ell_{s_1}^n)]_\theta
$$
by \cite[4.4.1]{BL}
$$ \norm{[\Phi_0^{n,m},\Phi_1^{n,m}]_\theta} \le 
\norm{\Phi_0^{n,m}}^{1-\theta}  \cdot  \norm{\Phi_1^{n,m}}^\theta. $$
Since by the interpolation theorem for $\ell_p(E)$'s 
together with the duality theorem \cite[4.5.2]{BL}
$$ [\ell_{s_0}^n(F_0'),\ell_{s_1}^n(F_1')]_\theta =
   \ell_{s_\theta}^n(F_\theta') \qquad \text{and} \qquad
[\ell_2^m(\ell_{s_0}^n),\ell_2^m(\ell_{s_1}^n)]_\theta =
\ell_2^m(\ell_{s_\theta}^n) $$
(isometrically), we have 
$$ \norm{\Phi_\theta^{n,m}} \le \norm{\LL(\ell_2^m, E_\theta) \hookrightarrow 
[\LL(\ell_2^m, E_0), \LL(\ell_2^m, E_1)]_\theta} \cdot 
\norm{[\Phi_0^{n,m},\Phi_1^{n,m}]_\theta}.$$
Consequently
\begin{align*}
 \mu_{s_\theta,2}(T:E_\theta \rightarrow F_\theta) &= 
\sup_{n,m} \norm{\Phi_\theta^{n,m}} \\ 
& \le\sup_{n,m} \{
   \norm{\LL(\ell_2^m, E_\theta) \hookrightarrow [\LL(\ell_2^m, E_0), 
\LL(\ell_2^m, E_1)]_\theta} \cdot \norm{\Phi_0^{n,m}}^{1-\theta} 
 \cdot  \norm{\Phi_1^{n,m}}^\theta \} \\ 
& \le d_\theta[E_0,E_1] \cdot \mu_{s_0,2}(T:E_0 \rightarrow F_0)^{1-\theta}
 \cdot \mu_{s_1,2}(T:E_1 \rightarrow F_1)^\theta, 
\end{align*}
the desired result. 
\qed
\\[10pt]
In the same way an analogous result for the $(r,2)$-summing norm can be 
obtained.
\par 
Applications of theorem~\ref{mixing} come from ``uniform estimates'' for
 $d_\theta[E_0,E_1]$. Pisier proved in \cite{pisier90} that
\begin{equation}
\label{koubal}
d_\theta[\ell_1,\ell_2]:= \sup_n d_\theta[\ell_1^n,\ell_2^n] 
< \infty;
\end{equation}
the proof is based on the Maurey factorization theorem which says that
 every operator $T:\ell_2 \rightarrow \ell_p$ ($1 \le p \le 2$) factorizes
 through an appropriate diagonal operator $D_\lambda: \ell_2 \rightarrow 
\ell_p$, and spaces of diagonal operators clearly behave well under 
interpolation. The fact \eqref{koubal} can also be obtained as an 
application (case (c) of the following estimates) of a deep result of Kouba 
\cite{kouba} on the complex interpolation of injective 
 tensor products of Banach spaces (see also 
\cite{dm99})---as a consequence of Kouba's work we 
get for a finite-dimensional interpolation couple $[E_0,E_1]$ and
 $0 \le \theta \le 1$ the following 
estimates for $d_\theta[E_0,E_1]$:
\begin{enumerate}[(a)]
\item
$d_\theta[E_0,E_1]=1$, if $E_0=E_1$.
\item 
$ d_\theta[E_0,E_1] \le \type (E_0')^{1-\theta} \cdot \type (E_1')^{\theta}.$ 
\item
$ d_\theta[E_0,E_1] \le (2/\sqrt{\pi}) \cdot \concave(E_0)^{1-\theta} \cdot
\concave(E_1)^{\theta},$
if the canonical bases in $E_0$ and $E_1$ are $1$-unconditional  
and induce the lattice structures.
\end{enumerate}
\vspace{-5pt}
Note that the constants given in (b) and (c) are different from
those in Kouba's work. They follow from a short analysis of his 
proofs in the finite-dimensional case:
 Since in our setting a
 Hilbert space is involved, Kouba's formulas 
(3.8) and (3.10) on p.~47--48 can both be
 changed into \mbox{$\gamma_z (T) \le \| \mspace{-1.5mu} | T \| 
 \mspace{-1.5mu} |_z$}. Moreover, calculating the terms $W^E(z)$ defined in
 Kouba's lemma~3.2 and lemma~3.3 (use two spaces instead of a family of 
Banach spaces, see also \cite[p.~218]{ccrsw}), 
one obtains 
$$W^E(z)= \type(E_0)^{1-\theta (z)} \cdot \type(E_1)^{\theta (z)} \text{ and }
W^E(z) = (2/\sqrt{\pi}) \cdot \con (E_0)^{1-\theta (z)} \cdot
\con (E_1)^{\theta (z)} $$
 (with $\theta (z)$ as in \cite[corollary~5.1]{ccrsw}), respectively.
This leads to the above estimates (note that $E_0$ and $E_1$ in Kouba's
 formulas, in our context have to be replaced by
 $E_0'$ and $E_1'$, and recall the 
well-known duality relation between $2$-concavity and $2$-convexity, see 
e.\,g. \cite[1.d.4]{lt}). 
\\[10pt] Later we also need a non-commutative version  of Pisier's result
 \eqref{koubal}.
Using an extension of Kouba's formulas for the Haagerup tensor product
 of operator spaces due to \cite{pisierop}, Junge in 
\cite[4.2.6]{junge} proved an analogue of \eqref{koubal} for Schatten classes:
\begin{equation}
\label{koubas}
d_\theta[\ui_1,\ui_2]:= \sup_n d_\theta[\ui_1^n,\ui_2^n] < \infty.
\end{equation}
\par
Finally we state a corollary on $\theta$-Hilbert spaces which together
 with \eqref{koubal} and \eqref{koubas} is crucial for our purposes.
\begin{cor}
\label{hilbert}
Let $0 \le \theta \le 1$, $E=[E_0,\ell_2^n]_\theta$ be a $n$-dimensional 
$\theta$-Hilbert space and \mbox{$2 \le s_\theta \le\infty$} given by 
$s_\theta =2/\theta$. Then 
$$ \mu_{s_\theta,2}(E \Id \ell_2^n) \le d_\theta[E_0,\ell_2^n] \cdot 
\pi_2(E_0 \Id \ell_2^n)^{1-\theta}.$$
\end{cor}
\section{\sc Bennett--Carl Inequalities for Symmetric Banach Sequence Spaces}
As announced the preceding interpolation theorem implies the Bennett--Carl
 result and its extension of Carl--Defant as an almost immediate consequence:
\begin{samepage}
\begin{cor}
\label{beca}
Let $1 \le u \le 2 $ and $1 \le u \le v \le \infty$. Then  for 
$2 \le s \le \infty $  such that $1/s= 1/2 - 1/u + \max(1/v,1/2)$
 $$ \sup_{n} \mu_{s,2}(\ell_u^n \Id \ell_v^n) < \infty.$$
In particular, for $2 \le r \le \infty$ such that $1/r = 1/u - \max(1/v,1/2)$
$$ \sup_{n} \pi_{r,2}(\ell_u^n \Id \ell_v^n) < \infty.$$
\end{cor}
\end{samepage}
\proof Only the case $1 \le u < v \le 2$ has to be considered; 
the case $2 \le v \le \infty$ then easily follows by factorization through
 $\ell_2^n$, and the case $u=v$ is trivial anyway.
 In what follows we use the complex interpolation formula for 
$\ell_p^n$'s without further mentioning. \\
i) Take first $v=2$. It is well-known (see e.\,g. \cite[22.4.8]{pietsch} 
or \eqref{hilfe1}) 
that
$$  \pi_2 (\ell_1^n \Id \ell_2^n)= 1. $$
For $1 \le u \le 2$ choose $0\le \theta \le 1$ such that 
$1/u = (1-\theta)/1 + \theta/2$. Then  $s_\theta:=2/\theta=u'$, and by
corollary~\ref{hilbert} together with \eqref{koubal}
$$ \mu_{u',2}(\ell_u^n \Id \ell_2^n) \le  d_\theta[\ell_1,\ell_2] 
< \infty.$$
ii) Let $1 \le u < v < 2$. Combining case~i),
$$ \mu_{u',2}(\ell_u^n \Id \ell_2^n) \le d_\theta[\ell_1,\ell_2],$$
and
$$ \mu_{2,2} (\ell_u^n \Id \ell_u^n ) = \norm{\ell_u^n \Id \ell_u^n } 
=1, $$ we have
$$ \mu_{s_{\tilde{\theta}},2}(\ell_u^n \Id \ell_v^n) \le 
\sup_n d_{\tilde{\theta}}[\ell_u^n,\ell_u^n] \cdot
d_\theta[\ell_1,\ell_2]^{1-\tilde{\theta}} < \infty,$$
with $\tilde{\theta} := (1/v -1/2)/(1/u - 1/2)$ and 
$1/{s_{\tilde{\theta}}} := (1-\tilde{\theta})/{u'} + \tilde{\theta}/2 
=1/2 - 1/u + 1/v=1/s$.
\qed
\par
As in the original proofs of Bennett and Carl the crucial step 
in the preceding proof is to show that
 for the symmetric Banach sequence space $E=\ell_u$
\begin{equation}
\label{crucial}
\sup_{n} \pi_{r,2} (E_n \Id \ell_2^n) < \infty,
\end{equation}
where $1 \le u \le 2$ and $1/r =1/u -1/2$. We now prove a result
within the framework of symmetric Banach sequence spaces which shows that 
\eqref{crucial} is sharp in a very strong sense.
 Take an arbitrary $2$-concave and $u$-convex
 Banach sequence space $E$---these geometric assumptions in particular 
imply that the continuous inclusions
 $\ell_u \subset E \subset \ell_2$ hold---which satisfies \eqref{crucial}. The 
following result shows that there is only one such space:
\begin{samepage}
\begin{theo}
\label{converse}
Let $1 \le u \le 2$ and $1/r =1/u -1/2$. For each $2$-concave and $u$-convex
 symmetric Banach sequence space $E$ the following are equivalent:
\em \vspace{-5pt}
\begin{enumerate}[(1)]
\item
$\sup_{n} \mu_{u',2}(E_n \Id \ell_2^n) <\infty$.
\item
$\sup_{n} \pi_{r,2}(E_n \Id \ell_2^n) < \infty$.
\item
$E=\ell_u$.
\end{enumerate}
\end{theo}
\end{samepage}
Clearly we only have to deal with the implication $(2) \Rightarrow (3)$; its
 proof is based on two lemmas. For the first one we invent the 
notion of ``enough symmetries in the orthogonal group''. 
Let $E=(\C^n,\norm{\cdot})$ be an $n$-dimensional Banach space. We say 
that $E$ has {\em enough symmetries in $\On(n)$} if there is a compact 
subgroup $G$ in $\On(n)$ such that
\begin{equation}
\label{invariant}
\forall \, u \in \LL(E) \, \forall \, g,g' \in G: \norm{u}=\norm{gug'}
\end{equation}
and
\begin{equation}
\label{enough}
\forall \, u \in \LL(E) \text{ with } ug=gu \text{ for all } g \in G
 \, \exists \, c \in \K: u = c \cdot \id_E.
\end{equation}
Basic examples of spaces with enough symmetries in the orthogonal group are 
the finite-dimensional spaces $E_n$ and $\ui_E^n$ associated to a symmetric
 Banach sequence space $E$.
 The following lemma extends 
the corresponding results in \cite[p.~233, 236]{CD97}.
\begin{lemma}
\label{2summing}
Let $E_n$ and $F_n$ have enough symmetries in $\On(n)$. Then 
\begin{equation}
\label{hilfe1}
\pi_2(E_n \Id F_n) = n^{1/2} \cdot 
\frac{\norm{\ell_2^n \Id F_n}}{\norm{\ell_2^n \Id E_n}},
\end{equation}
and for  $1 \le k \le n$
\begin{equation}
\label{hilfe2}
  \left ( \frac{n-k+1}{n} \right )^{1/2} \cdot
\frac{\norm{\ell_2^n \Id F_n}}{\norm{\ell_2^n \Id E_n}}
 \le x_k(E_n \Id F_n)
\le \left ( \frac{n}{k} \right )^{1/2} \cdot 
\frac{\norm{\ell_2^n \Id F_n}}{\norm{\ell_2^n \Id E_n}}.
\end{equation}
\end{lemma}
\proof \eqref{hilfe1}: Trace duality allows to deduce the lower estimate from 
the upper one:
$$
n \le \pi_2(\ell_2^n \Id F_n) \cdot \pi_2(F_n \Id \ell_2^n) 
  \le \norm{\ell_2^n \Id E_n} \cdot \pi_2(E_n \Id F_n) \cdot n^{1/2}
 \cdot \norm{\ell_2^n \Id F_n}^{-1}.
$$
For the proof of the upper estimate it may be assumed without loss of 
generality that $F_n=\ell_2^n$ (factorize through $\ell_2^n$). In this 
case it suffices to show that
$$
\norm{\ell_2^n \Id E_n}^{-1} \cdot B_{\ell_2^n}
$$
is John's ellipsoid $D_{\max}$ of maximal volume in $B_{E_n}$ (see e.\,g. 
\cite[3.8]{pisier89} or \cite[6.30]{djt}).
By definition there is a linear bijection $u: \ell_2^n \rightarrow E_n$ 
such that $u(B_{\ell_2^n}) = D_{\max}$. In particular, $\norm{u}=1$ and 
$N(u^{-1})=n$ ($N$ denotes the nuclear norm, see e.\,g. \cite[3.7]{pisier89}
 or \cite[6.30]{djt}). On the other hand by a standard averaging argument
 there is a linear bijection $v:\ell_2^n \rightarrow E_n$ with 
$\norm{v}=1$, 
$N(v^{-1})=n$ and $vg=gv$  for all $ g \in G$,
where $G$ is a compact group in $\On(n)$ satisfying \eqref{invariant} and
 \eqref{enough} (see \cite[3.5]{pisier89} which also holds in the complex
case). By property \eqref{enough} of
 $G$ and the fact that $\norm{v}=1$ we have
$v=\norm{\ell_2^n \Id E_n}^{-1} \cdot \id$. 
Then by Lewis' uniqueness theorem $v^{-1}u \in \On(n)$
 (\cite[3.7]{pisier89} or \cite[6.25]{djt}).
Altogether we finally obtain
$$
\norm{\ell_2^n \Id E_n}^{-1}\cdot B_{\ell_2^n}  =v(B_{\ell_2^n}) 
= v [v^{-1}u(B_{\ell_2^n})] 
= u(B_{\ell_2^n})=D_{\max}.
$$
\eqref{hilfe2}: Recall from \eqref{weyl} that 
$k^{1/2} \cdot x_k(T) \le \pi_2(T)$ for
 every $2$-summing operator $T$ acting between two Banach spaces. Together 
with \eqref{hilfe1} this gives the second inequality. The first then follows 
from the basic properties of the Weyl numbers (see e.\,g. \cite{koenig}):
\begin{align*}
1 &= x_n(\id_{\ell_2^n})  \\
  & \le x_k(\ell_2^n \Id F_n) \cdot x_{n-k+1}(F_n \Id \ell_2^n) \\
  & \le \norm{\ell_2^n \Id E_n} \cdot x_k(E_n \Id F_n) \cdot 
    \left (\frac{n}{n-k+1} \right )^{1/2} \cdot \norm{\ell_2^n \Id F_n}^{-1}. 
\end{align*}
\qed
\\[10pt]
The following obvious examples will be useful later. 
\begin{cor}
For $1 \le u,v \le \infty$ 
\begin{equation}
\label{formel}
\pi_2 (\ui_u^n \Id \ui_v^n) = 
n \cdot \frac{\max(1, n^{1/v -1/2})}{\max(1, n^{1/u -1/2})}
\end{equation}
and
\begin{equation}
\label{weylschatten}
x_{[n^2/2]} (\ui_u^n \Id \ui_v^n) \asymp 
\frac{\max(1, n^{1/v -1/2})}{\max(1, n^{1/u -1/2})}.
\end{equation}
\end{cor}
\par
 The preceding lemma turns out to be of special interest in combination 
with a result 
due to Szarek and Tomczak-Jaegermann \cite[proposition~2.2]{stj} which 
states that for each $2$-concave symmetric Banach sequence space $E$ 
\begin{equation}
\label{szarek}
 \norm{\ell_2^n \Id E_n} \asymp  n^{-1/2}
 \cdot \norm{\sum\nolimits_1^n e_i}_{E_n}.
\end{equation}
The second lemma which we need for the proof of theorem~\ref{converse}, 
is based on \eqref{szarek} and an important result 
about the interpolation of Banach lattices due to Pisier \cite{pisier79}.
\begin{lemma}
For $1 \le u \le 2$ let $E$ be a $u$-convex  and $u'$-concave symmetric 
Banach sequence space. Then
\begin{equation}
\label{zuerst}
\norm{E_n \Id \ell_u^n} \asymp \frac{n^{1/u}}{\norm{\sum_1^n e_i}_{E_n}}.
\end{equation}
In particular, if $E$ is even $2$-concave, then
\begin{equation}
\label{lu}
\norm{E_n \Id \ell_u^n} \asymp \frac{n^{1/u}}{\norm{\sum_1^n e_i}_{E_n}}
 \asymp \frac{n^{1/u-1/2}}{\norm{\ell_2^n \Id E_n}}.
\end{equation}
\end{lemma}
\proof \eqref{lu} follows directly from \eqref{zuerst} and \eqref{szarek},
 and clearly $n^{1/u} \le \norm{E_n \Id \ell_u^n} \cdot 
\norm{\sum_1^n e_i}_{E_n}$.
For the upper estimate in \eqref{zuerst} we only have to consider 
$1 < u <2$: The case
 $u=1$ is stated below in \eqref{border}, and a $2$-convex and $2$-concave
 symmetric Banach sequence space  necessarily equals $\ell_2$ with 
equivalent norms.
 Without loss
 of generality we may assume $\conu(E)= \mathbf{M_{(u')}}(E)=1$ (see
 \cite[1.d.8]{lt}). Then by \cite[theorem~2.2]{pisier79} there exists
 a symmetric Banach sequence space $E_0$ such that $E=[E_0,\ell_2]_\theta$
 with $\theta=2/{u'}$; moreover, we have $E_n=[E_0^n,\ell_2^n]_\theta$ with
 equal norms. The conclusion now follows by interpolation: It can be shown
 easily that 
\begin{equation}
\label{border}
\norm{E_0^n \Id \ell_1^ n} \le \frac{n}{\norm{\sum_1^n e_i}_{E_0^n}}
\end{equation}
 (see e.\,g. \cite[proposition~2.5]{stj}), hence
$$
\norm{E_n \Id \ell_u^n} \le \norm{E_0^n \Id \ell_1^n}^{1-\theta} 
\cdot \norm{\ell_2^n \Id \ell_2^n}^\theta \le 
\frac{n^{1-\theta}}{\norm{\sum_1^n e_i}_{E_0^n}^{1-\theta}}.
$$
Since $E_n=[E_0^n,\ell_2^n]_\theta$ is of $J$-type $\theta$ (i.\,e. 
$\norm{x}_{E_n} \le \norm{x}_{E_0^n}^{1-\theta} \cdot \norm{x}_{\ell_2^n}
^\theta$ for all $x \in E_n$), we have
$$
\norm{{\textstyle \sum_1^n e_i}}_{E_n} \le 
\norm{{\textstyle \sum_1^n e_i}}_{E_0^n}^{1-\theta} \cdot
 n^{\theta/2},
$$
and consequently
$$
\norm{E_n \Id \ell_u^n} \le \frac{n^{1-{\theta/2}}}{\norm{\sum_1^n e_i}_{E_n}}
 = \frac{n^{1/u}}{\norm{\sum_1^n e_i}_{E_n}}.
$$
\qed 
\\[10pt]
{\em Proof} of the implication $(2) \Rightarrow (3)$ in 
theorem~\ref{converse}:  
Assume that $\sup_{n} \pi_{r,2}(E_n \Id \ell_2^n) < \infty$. By 
\eqref{weyl}, \eqref{hilfe2} and \eqref{lu}
\begin{equation}
\label{lower}
\pi_{r,2}(E_n \Id \ell_2^n) \ge [n/2]^{1/r} \cdot x_{[n/2]}(E_n \Id \ell_2^n)
 \succ \frac{n^{1/r}}{\norm{\ell_2^n \Id E_n}} \asymp 
\norm{E_n \Id \ell_u^n},
\end{equation}
which by assumption shows that $\sup_n \norm{E_n \Id \ell_u^n} < \infty$.
 This clearly gives the claim.  \qed
\\[10pt] Note that \eqref{lower} does not depend on 
the special choice of $r$.
\par
If $E$ is a $2$-concave and $u$-convex $(1 \le u \le 2)$ symmetric
Banach sequence space
 different from $\ell_u$ (i.\,e. the inclusion $\ell_u \subset E$ is strict), 
then by theorem~\ref{converse} for $1/r = 1/u-1/2$ 
$$
\pi_{r,2}(E_n \Id \ell_2^n) \nearrow \infty.
$$
The following result gives the precise asymptotic order of the sequence
  \mbox{$(\pi_{r,2}(E_n \Id \ell_2^n))_n$}:
\begin{cor}
\label{precise}
For $1 \le u \le 2$ let $E$ be a $2$-concave and $u$-convex symmetric Banach
 sequence space. Then for $2 \le r,s \le \infty$ such that $1/r =1/u-1/2$ 
and $1/s=1/2-1/r$
$$
\pi_{r,2}(E_n \Id \ell_2^n) \asymp \mu_{s,2}(E_n \Id \ell_2^n) 
\asymp \frac{n^{1/r+1/2}}{\norm{\sum_1^n e_i}_{E_n}}.
$$ 
\end{cor}
\proof The lower estimate has already been shown in \eqref{lower}, and the
 upper estimate simply follows by factorization through $\ell_u^n$, the 
Bennett--Carl inequalities and \eqref{lu}. \qed
\begin{rem}
\label{remark}
\begin{enumerate}[(a)]
\item
Since a $u$-convex Banach lattice is $p$-convex for all $1\le p \le u$ 
(see \cite[1.d.5]{lt}), the formula in the preceding theorem even holds
 for all $2\le r \le \infty$ such that $1/r \ge 1/u -1/2$. 
\item 
For $1 \le u \le 2$ let $E$ be a $2$-concave and $u$-convex symmetric 
Banach sequence space, $F$ an arbitrary symmetric Banach sequence space,
and let $2 \le r \le \infty$ such that $1/r \ge
 1/u -1/2$. Then---by factorization through $\ell_2^n$ for
 the upper estimate and \eqref{lower} for the lower one---the following
 formula holds:
$$
\pi_{r,2}(E_n \Id F_n) \asymp n^{1/r} \cdot \frac{\norm{\ell_2^n \Id F_n}}
{\norm{\ell_2^n \Id E_n}};
$$
in particular, if $F$ is a $2$-concave, then
$$
\pi_{r,2}(E_n \Id F_n) \asymp n^{1/r} \cdot \frac{\norm{\sum_1^n e_i}_{F_n}}
{\norm{\sum_1^n e_i}_{E_n}}.
$$
Note that these results can be considered as extensions of \eqref{hilfe1}.
\item
For the special case $F=\ell_v$ ($1 \le u \le v \le 2$) the formulas in (b)
 even hold for all $2 \le r \le \infty$ such that $1/r \ge 1/u -1/v$; simply
 repeat the proof of corollary~\ref{precise} for $1/r =1/u -1/v$ 
and use the argument from remark (a).
\end{enumerate}
\end{rem}  
\section{\sc Bennett--Carl Inequalities for Unitary Ideals}
We now use Junge's counterpart \eqref{koubas} of \eqref{koubal} and our 
interpolation theorem~\ref{mixing} in order to show 
 a ``non-commutative'' analogue. Note first that for all 
$1 \le u,v \le \infty$ and $2 \le r \le \infty$
\begin{equation}
\label{low}
 n^{1/r} \le \pi_{r,2}(\ui_u^n \Id \ui_v^n),
\end{equation}
and hence also for $2 \le s \le \infty$
$$ n^{1/2 -1/s} \le \mu_{s,2}(\ui_u^n \Id \ui_v^n);$$
 this is a consequence of the
 trivial estimate $\pi_{r,2}(\ell_2^n \Id \ell_2^n) \ge n^{1/r}$ (insert
 $e_k$'s) and the fact that $\ell_2^n$ is $1$-complemented in each $\ui_u^n$
 (assign to each $x \in \ell_2^n$ the matrix $x \otimes e_1 \in \ui_u^n$).
 For $u,v$ considered in corollary~\ref{beca}
 this lower bound is optimal:  
\begin{cor}
\label{schatten}
Let $1 \le u \le 2 $ and $1 \le u \le v \le \infty$. Then for $2 \le s \le
 \infty$ such that $1/s= 1/2 - 1/u + \max(1/v,1/2)$
$$ \mu_{s,2}(\ui_u^n \Id \ui_v^n) \asymp n^{1/2-1/s}.$$
In particular, for $2 \le r \le \infty$ and $1/r = 1/u - \max(1/v,1/2)$
$$\pi_{r,2}(\ui_u^n \Id \ui_v^n) \asymp n^{1/r}.$$ 
\end{cor}
\proof  The proof of the upper bound is analogous to that of 
corollary~\ref{beca}:
Of course the complex interpolation formula for $\ui_p^n$'s is needed 
instead of that for $\ell_p^n$'s, and
in i) use $\pi_2(\ui_1^n \Id \ui_2^n) = n^{1/2}$ (see \eqref{formel}) and 
Junge's result \eqref{koubas} in order to obtain
$$
\mu_{u',2}(\ui_u^n \Id \ui_2^n) \le d_\theta[\ui_1^n,\ui_2^n] \cdot n^{(1-\theta)/2}
 \le d_\theta[\ui_1,\ui_2] \cdot n^{1/u-1/2},
$$
where $\theta=2/{u'}$.
Then in ii) one arrives at 
$$ \mu_{s_{\tilde{\theta}},2}(\ui_u^n \Id \ui_v^n) \prec  n^{(1-\tilde{\theta})(1/u - 1/2)}
 =n^{1/u-1/v}, $$
with $\tilde{\theta} := (1/v -1/2)/(1/u - 1/2)$ and 
$1/{s_{\tilde{\theta}}} = (1-\tilde{\theta})/{u'} + \theta/2 
=1/2 - 1/u + 1/v$=1/s.\qed
\par
Exploiting the ideas of the preceding section one easily obtains the 
asymptotic order of the $(r,2)$-summing and the $(s,2)$-mixing norm of
 identities between finite-dimensional unitary ideals $\ui_E^n$ and 
$\ui_2^n$:
\begin{cor}
\label{unitary}
For $1 \le u \le 2$ let $E$ be a $2$-concave and $u$-convex symmetric 
Banach sequence space. Then for all $2\le r,s \le \infty$ such that
 $1/r \ge 1/u -1/2$ and $1/s =1/2-1/r$
$$
\pi_{r,2}(\ui_E^n \Id \ui_2^n) \asymp 
\mu_{s,2}(\ui_E^n \Id \ui_2^n) \asymp \frac{n^{2/r+1/2}}{\norm{\sum_1^n 
e_i}_{E_n}}.
$$
\end{cor}
\proof Recall the simple fact that for all 
symmetric Banach sequence spaces $E$ and $F$
\begin{equation}
\label{ssui}
\norm{\ui_E^n \Id \ui_F^n} = \norm{E_n \Id F_n},
\end{equation}
and by the same reasoning as in remark~\ref{remark} (a) it is enough to deal 
with the case $1/r=1/u-1/2$. Then factorization through $\ui_u^n$ and 
\eqref{lu} give
$$
\mu_{u',2}(\ui_E^n \Id \ui_2^n) \prec \norm{\ui_E^n \Id \ui_u^n} \cdot n^{1/u-1/2}
 \asymp \frac{n^{2/u-1/2}}{\norm{\sum_1^n e_i}_{E_n}}
=\frac{n^{2/r+1/2}}{\norm{\sum_1^n e_i}_{E_n}},
$$
and in order to obtain the lower estimate apply again \eqref{hilfe2} together
 with \eqref{weyl} and the second asymptotic in \eqref{lu}:
$$
\pi_{r,2}(\ui_E^n \Id \ui_2^n) \ge [n^2/2]^{1/r} \cdot x_{[n^2/2]}(\ui_E^n
 \Id \ui_2^n) \succ \frac{n^{2/r}}{\norm{\ell_2^n \Id E_n}} \asymp 
\frac{n^{2/r+1/2}}{\norm{\sum_1^n e_i}_{E_n}}.
$$
\qed
\section{\sc Applications}
\subsection*{Weyl Numbers}
The results of the preceding sections can be used to improve
the estimates for Weyl numbers of identities on symmetric 
Banach sequence spaces and
 unitary ideals in \eqref{hilfe2}: The exponent $1/2$ in each of the two 
inequalities there can be replaced by $1/u-1/2$ whenever $u$-convexity
 and $2$-concavity assumptions are made.  
\begin{cor}
\label{weylcor}
For $1 \le u,v \le 2$ let $E$ and $F$ be  $2$-concave 
symmetric Banach sequence spaces where $E$ is $u$-convex and $F$ is $v$-convex.
Then there exist constants $C_u,C_v>0$ such that for all 
 $1 \le k \le n$ 
$$
C_v^{-1} \cdot \left (\frac{n-k+1}{n} \right )^{1/v-1/2} 
\cdot \frac{\norm{\sum_1^n e_i}_{F_n}}{\norm{\sum_1^n e_i}_{E_n}} \le
x_k(E_n \Id F_n) \le C_u \cdot \left (\frac{n}{k} \right )^{1/u-1/2} 
 \cdot \frac{\norm{\sum_1^n e_i}_{F_n}}{\norm{\sum_1^n e_i}_{E_n}},
$$
and all $1 \le k \le n^2$
$$
C_v^{-1} \cdot \left (\frac{n^2-k+1}{n^2} \right )^{1/v-1/2}
\cdot  \frac{\norm{\sum_1^n e_i}_{F_n}}{\norm{\sum_1^n e_i}_{E_n}}  \le 
x_k(\ui_E^n \Id \ui_F^n) 
\le C_u \cdot \left (\frac{n^2}{k} 
\right )^{1/u-1/2} 
 \cdot \frac{\norm{\sum_1^n e_i}_{F_n}}{\norm{\sum_1^n e_i}_{E_n}}.
$$
\end{cor}
\proof The upper estimates follow by using the inequality~\eqref{weyl} 
and the results from the preceding two sections, and the lower estimates 
then are immediate consequences of the upper ones---simply repeat the proof
 of \eqref{hilfe2} with a different exponent. \qed
\\[10pt]
Recall that for the embedding $\ell_u^n \Id \ell_2^n$, $1 \le u \le 2$ 
by \cite[2.3.3]{CD92} even the following equality is known:
 $x_k(\ell_u^n \Id \ell_2^n)=k^{1/2-1/u}$, $1 \le k \le n$. The second 
estimate in corollary~\ref{weylcor} implies that for $1 <u<2$
$$
x_k(\ui_u^n \Id \ui_2^n) \le C_u \cdot \left ( \frac{n}{k} \right)^{1/u-1/2},
$$
hence by \cite[2.3.2]{CD92}
$$
a_k(\ui_2^n \Id \ui_u^n) \ge C_u^{-1} \cdot \left 
(\frac{n^2-k+1}{n}\right)^{1/u-1/2}.
$$
This disproves the conjecture
$$
a_k(\ui_2^n \Id \ui_u^n) \asymp \max \left (1,\left 
(\frac{n^2-k+1}{n^2}\right)^{1/2} \cdot n^{1/u-1/2} \right) 
$$
from \cite[p.~249]{CD97} (put $k:=[n^2-n^\alpha+1]$, $1 < \alpha <2$).
\subsection*{Identities on Orlicz and Lorentz Sequence Spaces}
In the following we apply our results, in particular the 
corollaries~\ref{precise} and \ref{unitary}, 
 to two natural examples of 
symmetric sequence spaces: Orlicz and Lorentz sequence spaces 
(for their definition and basic properties we refer to \cite{lt77}). We only 
treat the case where the range space of the embedding is the 
finite-dimensional Hilbert space and leave the formulation for other spaces 
and the corollaries for Weyl numbers to the reader.
\par Let us start with 
Orlicz sequence spaces $\ell_M$. 
\begin{cor}
Let $1 <u < 2$ and $M$ be a strictly increasing Orlicz function which 
satisfies the $\Delta_2$-condition at zero. Assume that there exists 
$K>0$ such that for all $s,t \in (0,1]$
\begin{equation}
K^{-1} \cdot s^2 \le {M(st)}/{M(t)} \le K \cdot s^u.
\label{orlicz}
\end{equation} 
Then for $2 < r,s < \infty$ such that $1/r >1/u-1/2$ and $1/s=1/2 -1/r$
$$\pi_{r,2}(\ell_M^n \Id \ell_2^n) \asymp \mu_{s,2}(\ell_M^n \Id \ell_2^n)
 \asymp \frac{n^{1/r+1/2}}{\norm{\sum_1^n e_i}_{\ell_M^n}} 
\asymp n^{1/r+1/2} \cdot M^{-1}(1/n)
$$
and
$$\pi_{r,2}(\ui_{\ell_M}^n \Id \ui_2^n) \asymp \mu_{s,2}(\ui_{\ell_M}^n 
\Id \ui_2^n) 
\asymp \frac{n^{2/r+1/2}}{\norm{\sum_1^n e_i}_{\ell_M^n}}
\asymp n^{2/r+1/2} \cdot M^{-1}(1/n).
$$
\end{cor}
Note that \eqref{orlicz} together with the $\Delta_2$-condition assures that
 $\ell_M$ is $2$-concave and $p$-convex for all $1 \le p < u$ (see 
\cite[2.b.5]{lt}).
\par
Now we state an analogue for Lorentz sequence spaces $d(w,u)$. 
\begin{cor}
\label{lorentzcor}
Let $1 < u  <2$ and $w$ be such that 
$ n \cdot w_n^q  \asymp \sum_{i=1}^n w_i^q,$  where $q=2/(2-u)$. Then
for $2 < r,s < \infty$
 such that $1/r \ge 1/u -1/2$ and $1/s=1/2-1/r$
$$\pi_{r,2}(d_n(w,u) \Id \ell_2^n) \asymp \mu_{s,2}(d_n(w,u) \Id \ell_2^n)
 \asymp n^{1/r +1/2 -1/u} \cdot w_n^{-{1/u}}.
$$
and
$$\pi_{r,2}(\ui_{d(w,u)}^n \Id \ui_2^n) \asymp \mu_{s,2}(\ui_{d(w,u)}^n 
\Id \ui_2^n)
 \asymp n^{2/r +1/2 -1/u} \cdot w_n^{-{1/u}}.
$$
\end{cor}
Recall that the space $d(w,u)$ is $u$-convex, and if $1 \le u<2$, it is 
$2$-concave if and only if $w$ satisfies the condition in the assumption 
of the corollary 
(see \cite[p.~245--247]{reisner}).
\subsection*{Limit Orders}
Finally, we consider the asymptotic order of the sequences $(\pi_{r,2}
(\ui_u^n \Id \ui_v^n))_n$ for arbitrary $2 \le r \le \infty$, 
$1 \le u,v \le \infty$.
 Define the limit orders
 $$\lambda_\ell(\PP_{r,2},u,v) := \inf \{ \lambda>0 \, | \, 
\exists \, \rho>0 \, \forall \, n: \pi_{r,2}(\ell_u^n \Id \ell_v^n) \le \rho 
\cdot n^\lambda\}$$  
and 
$$\lambda_{\ui}(\PP_{r,2},u,v) := \inf \{ \lambda>0 \, | \, 
\exists \, \rho>0 \, \forall \,n: \pi_{r,2}(\ui_u^n \Id \ui_v^n) \le \rho \cdot 
n^\lambda\}.$$ 
Here we only handle the limit order of summing operators since---using the 
fact that $\PP_{r,2}$ and $\MM_{s,2}$ for $1/s + 1/r = 1/2$
 are almost equal---one can easily see that $\lambda_\ell(\PP_{r,2},u,v)=
 \lambda_\ell(\MM_{r,2},u,v)$ and
$\lambda_{\ui}(\PP_{r,2},u,v)= \lambda_{\ui}(\MM_{r,2},u,v)$
(with the obvious definition for the right 
sides of these equalities; see \cite[22.3.7]{pietsch}).
\par
The calculation of the limit order 
$\lambda_\ell(\PP_{r,2},u,v)$ was completed in \cite{CMP}:
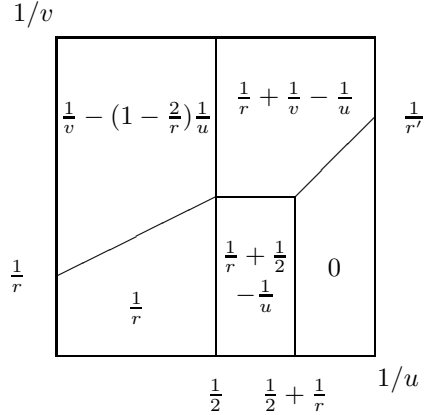
\begin{center}
\small
\unitlength3pt
\begin{picture}(60,50)
\put(10,10){\framebox(40,40)}
\put(30,10){\line(0,1){40}}
\put(30,30){\line(1,0){10}}
\put(40,10){\line(0,1){20}}
\put(40,30){\line(1,1){10}}
\put(10,20){\line(2,1){20}}
\put(40,10){\makebox(10,22){$0$}}
\put(30,10){\makebox(10,20){\shortstack{$\frac{1}{r}+ \frac{1}{2}$ \\
$-\frac{1}{u}$}}}
\put(30,35){\makebox(20,15){$\frac{1}{r} +
\frac{1}{v}-\frac{1}{u}$}}
\put(10,10){\makebox(20,12){$\frac{1}{r}$}}
\put(10,30){\makebox(20,20){$\frac{1}{v}-  (1-\frac{2}{r} )
 \frac{1}{u}$}}
\put(5,20){\makebox(0,0){$\frac{1}{r}$}}
\put(30,5){\makebox(0,0){$\frac{1}{2}$}}
\put(40,5){\makebox(0,0){$\frac{1}{2} + \frac{1}{r}$}}
\put(55,40){\makebox(0,0){$\frac{1}{r'}$}}
\put(7,53){\makebox(0,0){$1/v$}}
\put(53,7){\makebox(0,0){$1/u$}}
\end{picture}
\end{center}
Moreover, the proof in \cite{CMP} shows that the limit order is attained:
$\pi_{r,2}(\ell_u^n \Id \ell_v^n) \asymp n^{\lambda_\ell(\PP_{r,2},u,v)}$.
 In view of the results of section~4 the following conjecture seems to be
 natural:
\\[10pt]
{\bf Conjecture:} $\lambda_{\ui}(\PP_{r,2},u,v)=1/r +
 \lambda_\ell(\PP_{r,2},u,v)$.
\par 
 For the border cases $r=2$ (the $2$-summing norm) and $r=\infty$ 
(the operator norm) this conjecture  by \eqref{hilfe1} and 
\eqref{formel} is true.
In the following corollary we confirm the upper estimates of 
this conjecture for all $u,v$ and the lower ones for all $u,v$ except
 those in the upper left corner of the picture.
\begin{cor}
Let $1 \le u ,v\le $ and $2 < r < \infty$.
\vspace{-10pt}
\begin{enumerate}[(a)]
\item $\lambda_{\ui}(\PP_{r,2},u,v)=1/r + \lambda_\ell(\PP_{r,2},u,v)$ \quad 
for $1 \le u \le 2$.
\item 
$\lambda_{\ui}(\PP_{r,2},u,v) \le 1/r + \lambda_\ell(\PP_{r,2},u,v)$ \quad 
for $2 \le u \le \infty$,
with equality whenever \mbox{$1/v \le 1/r + (1- 2/r)(1/u)$.}
\end{enumerate}
\end{cor}
\proof Let $1/s:=1/2-1/r$. The upper estimates for the case 
$1 \le u \le 2$ follow from
 corollary~\ref{schatten}: Consider for $u_0 :=(1/2-1/r)^{-1}$ the 
following alternative: (i) $1/u \le 1/{u_0}$ or (ii) $1/u > 1/{u_0}$. Then
 the conclusion in case~(i) is a consequence of 
corollary~\ref{schatten} and the following factorization:
$$
\mu_{s,2}(\ui_u^n \Id \ui_v^n)  \le \norm{\ui_u^n \Id \ui_{u_0}^n} \cdot
 \mu_{s,2}(\ui_{u_0}^n \Id \ui_2^n) \cdot \norm{\ui_2^n \Id \ui_v^n} 
 \prec n^{2/r +1/2-1/u + \max(0,1/v-1/2)},
$$ 
and for (ii) look with $v_0:= (1/u-1/r)^{-1} \le 2$ at
\begin{align*}
\mu_{s,2}(\ui_u^n \Id \ui_v^n)  \le \mu_{s,2}(\ui_u^n \Id \ui_{v_0}^n) \cdot
 \norm{\ui_{v_0}^n \Id \ui_v^n} \prec n^{1/r + \max(0,1/r+1/v-1/u)}.
\end{align*}
Now let $2 \le u \le \infty$. Although this part is very 
close to the calculations made in \cite[Lemma~6]{CMP}, we give a short sketch 
 of the proof for the convenience of the reader. By \eqref{formel}
 and theorem~\ref{mixing} (with no interpolation in the range or the image),
\begin{align*}
\mu_{s,2}(\ui_u^n \Id \ui_2^n) \le \pi_2(\ui_u^n \Id \ui_2^n)^{2/r} \cdot 
\norm{\ui_u^n \Id \ui_2^n}^{1-{2/r}}= n^{1/r +1/2 -(1-2/r)(1/u)},
\end{align*}
hence, by factorization, for $1 \le v \le 2$
$$
\mu_{s,2}(\ui_u^n \Id \ui_v^n) \le n^{1/r + 1/v -(1-{2/r})(1/u)}.
$$
Furthermore, for $1/{v_1} := 1/r + (1-{2/r})(1/u)$
\begin{align*}
\mu_{s,2}(\ui_u^n \Id \ui_{v_1}^n )  \le \pi_2(\ui_u^n \Id \ui_2^n)^{2/r}
 \cdot \norm{\ui_u^n \Id \ui_u^n}^{1-{2/r}} = n^{2/r},
 \end{align*} 
hence
$$\mu_{s,2}(\ui_u^n \Id \ui_v^n) \le  n^{2/r} $$
for all $v_1 \le v \le \infty$. Finally, for all $2 < v < v_1$ and 
$0< \theta <1$ such that $1/v= (1-\theta)/{v_1} + \theta/2$
\begin{align*}
\mu_{s,2}(\ui_u^n \Id \ui_v^n) &\le \mu_{s,2}(\ui_u^n \Id \ui_{v_1}^n)
^{1-\theta} \cdot \mu_{s,2}(\ui_u^n \Id \ui_2^n)^\theta \\ 
& \le  n^{1/r + (1-\theta)/r + \theta(1/2 -(1-{2/r})(1/u))}  
 =  n^{1/r + 1/v -(1-{2/r})(1/u)}. 
\end{align*}
Looking at the picture for $\lambda_\ell(\Pi_{r,2}, u,v)$ one can see that
 these are the desired results.
For the lower estimates recall \eqref{weyl}:
$$ [{n^2}/2 ]^{1/r} \cdot x_{[{n^2}/2 ]}(\ui_u^n \Id \ui_v^n) \le 
\pi_{r,2}(\ui_u^n \Id \ui_v^n),$$
hence \eqref{weylschatten} implies
$$\pi_{r,2}(\ui_u^n \Id \ui_v^n) \succ 
\begin{cases}
n^{2/r + 1/v - 1/u} & \text{if $1 \le u,v \le 2$,}\\
n^{2/r + 1/2 -1/u} & \text{if $1 \le u \le 2 \le v \le \infty$,} \\
n^{2/r} & \text{if $2 \le u,v \le \infty$.} 
\end{cases}
$$
Using \eqref{low}, these estimates can be improved for those 
$u,v$ for which $\lambda_\ell(\PP_{r,2},u,v)=0$.
\qed 
\par
Our results for $\lambda_{\ui}(\PP_{r,2},u,v)$ can be summarized in the
following picture:
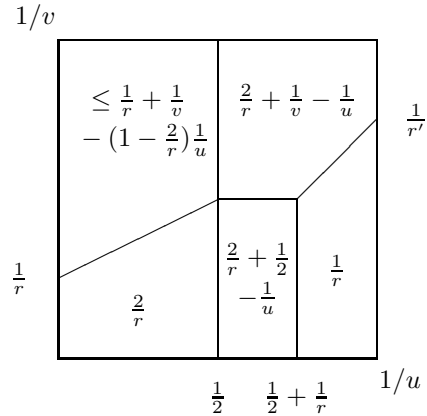
\begin{center}
\small
\unitlength3pt
\begin{picture}(60,60)
\put(10,10){\framebox(40,40)}
\put(30,10){\line(0,1){40}}
\put(30,30){\line(1,0){10}}
\put(40,10){\line(0,1){20}}
\put(40,30){\line(1,1){10}}
\put(10,20){\line(2,1){20}}
\put(40,10){\makebox(10,22){$\frac{1}{r}$}}
\put(30,10){\makebox(10,20){\shortstack{$\frac{2}{r}+ \frac{1}{2}$ \\
$-\frac{1}{u}$}}}
\put(30,35){\makebox(20,15){$\frac{2}{r} +
\frac{1}{v}-\frac{1}{u}$}}
\put(10,10){\makebox(20,12){$\frac{2}{r}$}}
\put(10,30){\makebox(20,20){\shortstack{$ \le \frac{1}{r} + \frac{1}{v}$ \\
$ \text{ } - (1-\frac{2}{r}) \frac{1}{u}$}}}
\put(5,20){\makebox(0,0){$\frac{1}{r}$}}
\put(30,5){\makebox(0,0){$\frac{1}{2}$}}
\put(40,5){\makebox(0,0){$\frac{1}{2} + \frac{1}{r}$}}
\put(55,40){\makebox(0,0){$\frac{1}{r'}$}}
\put(7,53){\makebox(0,0){$1/v$}}
\put(53,7){\makebox(0,0){$1/u$}}
\end{picture}
\end{center}
\par
The contents of this article will be part of the second named author's 
thesis written at the Carl von Ossietzky University of Oldenburg under
 the supervision of the first named author.
%

\begin{thebibliography}{CCRSW82}

\bibitem[Ben73]{bennett}
G.~Bennett, \emph{Inclusion mappings between $\ell^p$-spaces}, J. Funct. Anal.
  \textbf{13} (1973), 20--27.

\bibitem[BL78]{BL}
J.~Bergh and J.~L\"ofstr\"om, \emph{Interpolation spaces}, Springer-Verlag,
  1978.

\bibitem[Car74]{carl}
B.~Carl, \emph{Absolut $(p,1)$-summierende identische {O}peratoren von $\ell_u$
  nach $\ell_v$}, Math. Nachr. \textbf{63} (1974), 353--360.

\bibitem[CD92]{CD92}
B.~Carl and A.~Defant, \emph{Tensor products and {G}rothendieck type
  inequalities of operators in {$L_p$}-spaces}, Trans. Amer. Math. Soc.
  \textbf{331} (1992), 55--76.

\bibitem[CD97]{CD97}
B.~Carl and A.~Defant, \emph{Asymptotic estimates for approximation quantities
  of tensor product identities}, J. Approx. Theory \textbf{88} (1997),
  228--256.

\bibitem[CMP78]{CMP}
B.~Carl, B.~Maurey, and J.~Puhl, \emph{{G}renzordnungen von
  absolut-$(r,p)$-summierenden {O}peratoren}, Math. Nachr. \textbf{82} (1978),
  205--218.

\bibitem[CCRSW82]{ccrsw}
R.~R. Coifman, M.~Cwikel, R.~Rochberg, Y.~Sagher, and G.~Weiss, \emph{A theory
  of complex interpolation for families of {B}anach spaces}, Advances in Math.
  \textbf{43} (1982), 203--229.

\bibitem[DF93]{df}
A.~Defant and K.~Floret, \emph{Tensor norms and operator ideals},
  North-Holland, 1993.

\bibitem[DM99]{dm99}
A.~Defant and C.~Michels, \emph{A complex interpolation formula for 
tensor products of vector-valued {B}anach function spaces}, submitted (1999).

\bibitem[DJT95]{djt}
J.~Diestel, H.~Jarchow, and A.~Tonge, \emph{Absolutely summing operators},
  Cambridge Studies in Advanced Mathematics 43, 1995.

\bibitem[Jun96]{junge}
M.~Junge, \emph{Factorization theory for spaces of operators}, Univ.
  Kiel, Habilitationsschrift, 1996. Currently available on the net under \\
{\tt 
http://www-computerlabor.math.uni-kiel.de/\symbol{126}mjunge/preprints.html}

\bibitem[K{\"o}n86]{koenig}
Her. K{\"o}nig, \emph{Eigenvalue distributions of compact operators},
  Birkh\"auser, 1986.

\bibitem[Kou91]{kouba}
O.~Kouba, \emph{On the interpolation of injective or projective tensor products
  of {B}anach spaces}, J. Funct. Anal. \textbf{96} (1991), 38--61.

\bibitem[LT77]{lt77}
J.~Lindenstrauss and L.~Tzafriri, \emph{Classical {B}anach spaces {I}: Sequence
  spaces}, Springer-Verlag, 1977.

\bibitem[LT79]{lt}
J.~Lindenstrauss and L.~Tzafriri, \emph{Classical {B}anach spaces {II}:
  Function spaces}, Springer-Verlag, 1979.

\bibitem[MM99]{mm99}
L.~Maligranda and M.~Masty\l o, \emph{Inclusion mappings between
 Orlicz sequence spaces}, preprint (1999).

\bibitem[Mau74]{maurey}
B.~Maurey, \emph{Th\'eor\`emes de factorisation pour les op\'erateurs
  lin\'eaires \`a valeurs dans les espaces {$L^p$}}, Ast\'erisque 11, 1974.

\bibitem[Pie80]{pietsch}
A.~Pietsch, \emph{Operator ideals}, North-Holland, 1980.

\bibitem[Pie87]{pietsch2}
A.~Pietsch, \emph{Eigenvalues and $s$-numbers}, Cambridge Studies in Advanced
  Mathematics 13, 1987.

\bibitem[PT68]{pt}
A.~Pietsch and H.~Triebel, \emph{Interpolationstheorie f\"ur {B}anachideale von
  beschr\"ankten linearen {O}peratoren}, Studia Math. \textbf{31} (1968),
  95--109.

\bibitem[Pi79]{pisier79}
G.~Pisier, \emph{Some applications of the complex interpolation method to
  {B}anach lattices}, J. Analyse Math. \textbf{35} (1979), 264--281.

\bibitem[Pi89]{pisier89}
G.~Pisier, \emph{The volume of convex bodies and {B}anach space geometry},
  Cambridge Tracts in Math. 94, 1989.

\bibitem[Pi90]{pisier90}
G.~Pisier, \emph{A remark on {$\Pi_2(\ell_p,\ell_p)$}}, Math. Nachr. 
\textbf{148} (1990), 243--245.

\bibitem[Pi96]{pisierop}
G.~Pisier, \emph{The operator {H}ilbert space {$OH$}, complex interpolation and
  tensor norms}, Mem. Amer. Math. Soc. 585, 1996.

\bibitem[Rei81]{reisner}
S.~Reisner, \emph{A factorization theorem in {B}anach lattices and its
  applications to {L}orentz spaces}, Ann. Inst. Fourier \textbf{31} (1981),
  239--255.

\bibitem[STJ80]{stj}
S.~Szarek and N.~Tomczak-Jaegermann, \emph{On nearly {E}uclidean decomposition
  for some classes of {B}anach spaces}, Compositio Math. \textbf{40} (1980),
  367--385.

\bibitem[TJ89]{tj}
N.~Tomczak-Jaegermann, \emph{Banach--{M}azur distances and finite-dimensional
  operator ideals}, Longman Scientific \& Technical, 1989.

\end{thebibliography}
\newcommand{\etalchar}[1]{$^{#1}$}
\providecommand{\bysame}{\leavevmode\hbox to3em{\hrulefill}\thinspace}

\end{document}